\documentclass[12pt]{amsart}

\usepackage{amsfonts}
\usepackage{mathrsfs}
\usepackage{amssymb}

\newcommand{\bone}{\, \mathbf{1}}
\newcommand{\bD}{\mathbb{D}}
\newcommand{\bE}{\mathbb{E}}
\newcommand{\bP}{\mathbb{P}}
\newcommand{\bK}{\mathbb{K}}
\newcommand{\bR}{\mathbb{R}}
\newcommand{\bQ}{\mathbb{Q}}
\newcommand{\bZ}{\mathbb{Z}}
\newcommand{\cF}{\mathcal{F}}
\newcommand{\cG}{\mathcal{G}}
\newcommand{\cH}{\mathcal{H}}

\newcommand{\epX}{\varepsilon(X)}

\newcommand{\Verty}{\Vert _\infty}
\newcommand{\VertG}{\Vert _\cG}
\newcommand{\eps}{\varepsilon}
\newcommand{\convinfty}{\stackrel{L^\infty}{\longrightarrow}}
\newcommand{\ess}{\mathrm{ess}}
\newcommand{\supp}{\mathrm{supp}}

\theoremstyle{plain} \newtheorem{Theo}{Theorem}[section]
\theoremstyle{plain} \newtheorem{Lemma}[Theo]{Lemma}
\theoremstyle{plain} \newtheorem{Cor}[Theo]{Corollary}
\theoremstyle{plain} \newtheorem{Prop}[Theo]{Proposition}
\theoremstyle{remark} \newtheorem{Rem}[Theo]{Remark}
\theoremstyle{definition} \newtheorem{Exm}[Theo]{Example}
\theoremstyle{definition} \newtheorem{Def}[Theo]{Definition}
\theoremstyle{remark} \newtheorem{Not}[Theo]{Notation}

\title[Local field expectation]{Expectation, Conditional Expectation 
and Martingales in Local Fields}
\author{Steven N.\ Evans}
\address{
  Department of Statistics \#3860 \\
  University of California at Berkeley \\
  367 Evans Hall \\
  Berkeley, CA 94720-3860 \\
  U.S.A.}
\email{evans@stat.berkeley.edu}
\urladdr{http://www.stat.berkeley.edu/users/evans/}
\thanks{SNE supported in part by NSF grant DMS-0405778.}
\author{Tye Lidman}
\email{tlid@berkeley.edu}

\begin{document}

\begin{abstract}
We investigate a possible definition of expectation 
and conditional expectation for random variables with values
in a local field such as the $p$-adic numbers.  
We define the expectation by analogy with the observation
that for real-valued random variables in $L^2$ the expected value is the 
orthogonal projection onto the constants.
Previous work has shown that the local field version of
$L^\infty$ is the appropriate counterpart of
$L^2$, and so the expected value of 
a local field-valued random variable is defined to
be its ``projection'' in $L^\infty$ onto the constants.  Unlike the real
case, the resulting projection is not typically a single constant, but rather
a ball in the metric on the local field.  However, many properties of this
expectation operation and the corresponding conditional expectation mirror
those familiar from the real-valued case; for example, conditional
expectation is, in a suitable sense, a contraction on $L^\infty$ and
the tower property holds.
We also define the corresponding notion of martingale, show that several
standard examples of martingales (for example, sums or products of suitable
independent random variables or ``harmonic'' functions composed with Markov chains)
have local field analogues, and obtain versions of the optional sampling
and martingale convergence theorems.
\end{abstract}
\maketitle

\section{Introduction}

Expectation and conditional expectation of real-valued random variables
(or, more generally, Banach space-valued random variables) and 
the corresponding notion of martingale are fundamental objects of probability theory.
In this paper we investigate whether there are analogous notions
for random variables with values in a local field 
(that is, a locally compact, non-discrete, totally 
disconnected, topological field) -- 
a setting that shares the linear structure which underpins many of
the properties of the classical entities.

The best known example of a local field
is the field of {\em $p$-adic numbers}
for some positive prime $p$.  This field is defined as follows.
We can write any non-zero rational number $r \in \bQ \backslash \{ 0 
\}$ uniquely as $r=p^s (a/b)$, with $a,b,$ and $s$ integers,
where $a$ and $b$ are not divisible by $p$. Set $| r | = p^{-s}$.
If we set $| 0 | =0$, then the map $| \cdot |$ has
the properties:
\begin{equation}
\label{valuation}
\begin{split}
| x | & = 0 \Leftrightarrow x=0 \\
| xy | & = | x | | y | \\
| x+y | & \le | x | \vee | y |.
\end{split}
\end{equation}
The map $(x,y) \mapsto | x-y |$ defines a metric on $\bQ$ and we
denote the completion of $\bQ$ in this metric by $\bQ_p$.
The field operations on $\bQ$ extend continuously  to make $\bQ_p$
a  topological field called the {\em $p$-adic numbers}.
The map $| \cdot |$ also extends continuously and the extension
has properties (\ref{valuation}).  

The closed unit ball around $0$,
$\bZ_p = \{x \in \bQ_p : |x| \le 1 \}$,
is the closure in $\bQ_p$ of the integers $\bZ$,
and is thus
a ring (this is also apparent from (\ref{valuation})),
called the  {\em $p$-adic integers}.
As
$\bZ_p = \{x \in \bQ_p : |x| < p\}$, the set $\bZ_p$ is also open.
Any other ball around $0$ is of the form
$\{x \in \bQ_p : |x| \le p^{-k}\} = p^k \bZ_p$
for some integer $k$.  

Every local field is either a finite algebraic extension of the 
$p$-adic number field for some prime $p$ or a finite algebraic 
extension of the $p$-{\em series field}; that is, the field of formal 
Laurent series with coefficients drawn from the finite field with $p$ 
elements.)  A  locally compact, non-discrete, 
topological field that
is not totally disconnected is necessarily
either the real or the complex numbers.

   From now on, we let $\bK$ be a fixed local field.  Good general
reference for the properties of local fields and analysis on them are
\cite{MR791759, MR512894, MR0487295}.  The following are the properties we need.

There is a real-valued mapping $x \mapsto | x |$ on $\bK$ 
called the non-archimedean valuation
with the properties (\ref{valuation}).
The third of these properties is  the {\em ultrametric inequality} or
the {\em strong triangle inequality}.
The map $(x,y) \mapsto | x-y |$ on $\bK \times \bK$ is a metric on $\bK$
which gives the topology of $\bK$.
A consequence of of the strong triangle inequality is that if
$| x | \ne | y |$, then $| x+y | = | x | \vee | y |$.
This latter result implies that for every ``triangle''
$\{ x,y,z \} \subset \bK$ we have that at least two of the lengths $| x-y |$,
$| x-z |$, $| y-z |$ must be equal and is therefore often
called the {\em isosceles triangle property}.

The valuation takes the values
$\{ q^k : k \in \bZ \}
\cup \{ 0 \}$, where $q = p^c$ for some prime $p$ and positive integer $c$
(so that for $\bK=\bQ_p$ we have $c=1$).
Write $\bD$ for $\{ x \in \bK : | x | \le 1 \}$ (so that $\bD=\bZ_p$
when $\bK=\bQ_p$).  
Fix $\rho \in \bK$ so that $| \rho | =q^{-1}$.  Then
\begin{equation*}
\rho^k \bD  = \{ x: | x | \le q^{-k} \} = \{ x: | x | < q^{-(k-1)} \}
\end{equation*}
for each $k \in \bZ$ (so that for $\bK=\bQ_p$ we could take $\rho = p$).
The set $\bD$ is the unique maximal compact subring  of $\bK$
(the  {\em ring of integers} of 
$\bK$). Every ball in $\bK$
is of the form  $x + \rho^k \bD$ for some $x \in \bD$ and $k \in \bZ$.
If $B = x + \rho^k \bD$ and $C = y + \rho^\ell \bD$ are two such balls,
then 
\begin{itemize}
\item
$B \cap C = \emptyset$, if $|x-y| > q^{-k} \vee q^{-\ell}$,
\item
$B \subseteq C$, if $|x-y| \vee q^{-k} \le q^{-\ell}$,
\item
$C \subseteq B$, if $|x-y| \vee q^{-\ell} \le q^{-k}$.
\end{itemize}
In particular, if $q^{-k} = q^{-\ell}$, then either $B \cap C = \emptyset$
or $B=C$, depending on whether or not $|x-y| > q^{-k} = q^{-\ell}$
or $|x-y| \le q^{-k} = q^{-\ell}$.

We have shown in a sequence papers \cite{MR990478, MR1118442,
MR1245397, MR1414928, MR1832433, MR1873668, MR1934156, math.PR/0602478}
that the natural analogues on $\bK$
of the centered Gaussian measures on $\bR$
are the normalized restrictions of Haar measure
on the additive group of $\bK$ to the compact
the balls $\rho^k \bD$ and the point mass at $0$.  
There is a significant literature on probability on the
$p$-adics and other local fields.  The above papers
contain numerous references to this work, much of which concerns
Markov processes taking values in local fields.  There are
also extensive surveys of the literature in the books
\cite{MR1746953, MR1848777, MR2105195}.

It is not immediately clear how one should approach defining the expectation
of a local field valued random variable $X$.  Even if $X$ only
takes a finite number of values $\{x_1, x_2, \ldots, x_n\}$, then
the object $\sum_k x_k \bP\{X=x_k\}$ doesn't make any sense because
$x_k \in \bK$ whereas $\bP\{X=x_k\} \in \bR$.  However, it is
an elementary fact that if $T$ is a real-valued random variable
with $\bE[T^2]< \infty$, then $c \mapsto \bE[(T-c)^2]$ is uniquely
minimized by $c = \bE[T]$.  Of course, since this
observation already uses the notion of expectation it does not lead to
an alternative way of defining the expected value of a real-valued
random variable.  Fortunately, we can do something similar, but non-circular,
in the local field case.

Fix a probability space $(\Omega, \cF, \bP)$.  By a $\bK$-valued
random variable, we mean a measurable map from $\Omega$ equipped
with $\cF$ into $\bK$ equipped with its Borel $\sigma$-field.
Let $L^\infty$ be the space of
$\bK$-valued random variables $X$ that satisfy
$\|X\|_\infty := \ess \sup |X| < \infty$.  It is clear that 
$L^\infty$ is a vector space over $\bK$.  If we identify
two random variables as being equal when they are equal almost surely, then
\[
\begin{split}
\| X \|_\infty & = 0 \Leftrightarrow X=0 \\
\| c X \|_\infty & = | c | \| X \|_\infty, \quad c \in \bK, \\
\| X+Y \|_\infty & \le \| X \|_\infty \vee \| Y \|_\infty.
\end{split}
\]
The map $(X,Y) \mapsto \|X-Y\|_\infty$ defines a metric on $L^\infty$
(or, more correctly, on equivalence classes under the relation
of equality almost everywhere), and $L^\infty$ is complete in this metric.  
Hence $L^\infty$ is an instance of a Banach algebra over $\bK$.
 
 It is apparent from the papers on analogues of Gaussian measures
cited above that $L^\infty$ is the natural local field counterpart of the
real Hilbert space $L^2$.  In particular, there is a natural notion
of orthogonality on $L^\infty$ (albeit one which does not come from an
inner product structure).

\begin{Def}
Given $X \in L^\infty$, set
$\varepsilon(X) = \inf\{\Vert X-c \Verty : c \in \bK\}$. 
The {\em expectation} of the $\bK$-valued random variable 
$X$ is the subset of $\bK$ given by
\[
\bE[X]
:=
\{c \in \bK : \Vert X-c \Verty = \varepsilon(X)\}
=
\{c \in \bK : \Vert X-c \Verty \leq \varepsilon(X)\}.  
\]
\end{Def}

We show in Section \ref{S:expectation_exists} that $\bE[X]$ is non-empty.
Note that if $c' \in \bE[X]$ and $c'' \in \bK$ is such that
$|c''-c'| \le \varepsilon(X)$, then,
by the strong triangle inequality, $c'' \in \bE[X]$. Thus $\bE[X]$
is a (closed) ball in $\bK$ (where we take a single point as being a ball).

Observe that we use the same notation for expectation of
$\bK$-valued and $\bR$-valued random variables.  This
should cause no confusion: we  either indicate explicitly whether
a random variable has values in $\bK$ or $\bR$, or this
will be clear from context.

The outline of the rest of the paper is the following.
We show in Section \ref{S:expectation_exists} that the
expected value of a random variable in $L^\infty$ is
non-empty, remark on some of the properties of the expectation
operator, and motivate the definition of conditional expectation
by considering the situation where the conditioning $\sigma$-field
is finitely generated or, more generally,
has an associated regular conditional probability. The appropriate
definition of the conditional expectation of $X \in L^\infty$
given a sub-$\sigma$-field $\cG \subseteq \cF$ is not, as one might first
imagine, the $L^\infty$ projection of $X$ onto $L^\infty(\cG)$ ($:=$ the
subspace of $L^\infty$ consisting of $\cG$-measurable random variables).
For this reason, we need to do some preparatory work in
Sections \ref{S:cond_ess_sup} and \ref{S:cond_Linfty} before finally
presenting the construction of conditional expectation in 
Section \ref{S:constr_cond_exp} and describing its elementary properties
in Section \ref{S:prop_cond_exp}.  We establish an analogue of the
``tower property'' in Section \ref{S:tower_prop} and obtain a
counterpart of the fact for classical conditional
expectation that conditioning is a contraction on $L^2$
(both of these results need to be suitably interpreted due
to the conditional expectation being typically a set of random
variables rather than a single one).  We introduce the associated
notion of martingale in Section \ref{S:mart} and observe that
several of the classical examples of martingales have local field
analogues.  We develop counterparts of the optional sampling theorem
and martingale convergence theorem in Sections \ref{S:opt_sampling}
and \ref{S:mart_conv}, respectively.

\medskip
\noindent{\bf Note:} We adopt the convention that all equalities and
inequalities between random variables should be interpreted as holding
$\bP$-almost surely.

\section{Expectation}
\label{S:expectation_exists}

	\begin{Theo}
	\label{T:expectation_exists}
	The expectation of a random variable $X \in L^\infty$ is non-empty.
It is the smallest closed ball in $\bK$  that contains $\supp{X}$
	(the closed support of $X$).
	\end{Theo}

	\begin{proof}
By the strong triangle inequality
$\Vert X-c \Verty \le \Vert X \Verty \vee |c|$, and
$\Vert X-c \Verty =  |c|$ for $|c| > \Vert X \Verty$.
Therefore, the infimum of $c \mapsto \Vert X-c \Verty$
over all $c \in \bK$ is the same as the infimum over 
$\{c \in \bK : |c| \le \Vert X \Verty\}$ and any point
$c \in \bK$ at which the infimum of
is achieved must necessarily
satisfy $|c| \le \Vert X \Verty$.   That is,
$\varepsilon(X) = \inf\{\Vert X-c \Verty : |c| \le \Vert X \Verty\}$
and 
$\bE[X]
=
\{c: |c| \le  \Vert X \Verty, \;  \Vert X-c \Verty = \varepsilon(X)\}.
$

Again by the strong triangle inequality, 
the function $c \mapsto \Vert X-c \Verty$
is continuous.
Consequently,
$\bE[X]$
is non-empty as the set of points at which a continuous function
on a compact set attains its infimum.  

As we observed in the Introduction, $\bE[X]$ is a ball of radius 
($=$ diameter) $\varepsilon(X)$.  
If $x \in \supp X$ is not in $\bE[X]$ and $c$ is any point in $\bE[X]$,
then, by the strong triangle inequality, $|x-c| > \varepsilon(X)$
and $\Vert X-c \Verty > \varepsilon(X)$, contradicting the definition
of $\bE[X]$.  Thus $\supp X \subseteq \bE[X]$.  
Hence, if the smallest ball containing $\supp X$ is not $\bE[X]$, it
must be a ball contained in $\bE[X]$ with diameter 
$r < \varepsilon(X)$. However, if $c$ is any point contained in 
the smaller ball, then $|x-c| \le r$ for all $x \in \supp X$, contradicting
the definition of $\varepsilon(X)$.
\end{proof}

Our notion of
expectation shares some of the features of both
the mean and the variance of a real-valued variable.
Any point in the ball $\bE[X]$ is as good a single summary
of the ``location'' of $X$ as any other, whereas the
diameter of $\bE[X]$ (that is, $\varepsilon(X)$) is
a measure of the ``spread'' of $X$.

Some properties of $\bE[X]$ are immediate.
It is easily seen that for constants 
$k,b \in \bK$, $\bE[k X+b]=k\bE[X]+b$.  We do not have complete linearity, however, since $\bE[X+Y]$ is only a subset of $\bE[X]+\bE[Y]$, with equality when $X$ and $Y$ are independent.
This follows from
the fact that $\supp (X+Y) \subseteq \supp X + \supp Y$, with equality
when $X$ and $Y$ are independent.  
Also, if $X$ and $Y$ are independent, then $\bE[XY]=\bE[X]\bE[Y]$.
These remarks further support
our assertion that $\bE[X]$ combines the properties of the
mean and the variance for real-valued random variables.

Define the Hausdorff distance between two subsets $A$ and $B$
of $\bK$ to be
\[
d_H(A,B) :=
\sup_{a \in A} \inf_{b \in B}|a-b| 
\vee
\sup_{b \in B} \inf_{a \in A} |b-a|.
\]
We know from Theorem \ref{T:expectation_exists}
that $\bE[X]$ and $\bE[Y]$ are balls with diameters
$\varepsilon(X)$ and $\varepsilon(Y)$, respectively. We have
one of the alternatives $\bE[X] = \bE[Y]$, 
$\bE[X] \subsetneq \bE[Y]$, $\bE[Y] \subsetneq \bE[X]$,
or $\bE[X] \cap \bE[Y] = \emptyset$.  Suppose that
$\bE[X] \subsetneq \bE[Y]$, so that $\supp X \subseteq \bE[X]$
and there exists $y \in \supp Y$ such that $y$ is not in the
unique ball of diameter $q^{-1} \varepsilon(Y)$ containing $\bE[X]$.
Then, by the strong triangle inequality,
$|x-y| = \varepsilon(Y)$ for all $x \in \supp X$, and so
$d_H(\supp X, \supp Y) \ge \varepsilon(Y) = d_H(\bE[X], \bE[Y])$
in this case.  Similar arguments in the other cases show that
\[
d_H(\bE[X], \bE[Y]) \le d_H(\supp X, \supp Y) \le \Vert X - Y \Verty.
\]
This is analogous to the continuity of real-valued expectation with
respect to the real $L^p$ norms.

Rather than develop more properties of expectation, we move
on to the corresponding definition of conditional expectation
because, just as in the real case, expectation is the
special case of conditional expectation that occurs when the
conditioning $\sigma$-field is the trivial $\sigma$-field
$\{\emptyset, \Omega\}$, and so results for expectation are
just special cases of ones for conditional expectation.

In order to motivate the definition of conditional expectation, first
consider the special case when the conditioning 
$\sigma$-field $\cG \subseteq \cF$ is generated by a finite partition 
$\{A_1, A_2, \ldots, A_n\}$ of $\Omega$.  In line
with our definition of $\bE[X]$, a reasonable definition
of $\bE[X \, | \, \cG]$ would be the set of 
$\cG$-measurable random variables $Y$ such that for each $k$
the common value of $c_k := Y(\omega)$ for $\omega \in A_k$ satisfies
\[
\ess \sup \{|X(\omega) - c_k| : \omega \in A_k\}
=
\inf_{c \in \bK} 
\ess \sup \{|X(\omega) - c| : \omega \in A_k\}.
\]
Equivalently, suppose we define $\varepsilon(X,\cG)$
to be the $\cG$-measurable, $\bR$-valued random variable that
takes the value 
$\inf_{c \in \bK} 
\ess \sup \{|X(\omega) - c| : \omega \in A_k\}$ on $A_k$, then
$\bE[X \, | \, \cG]$ is the set of 
$\cG$-measurable random variables $Y$ such that
$|X-Y| \le \varepsilon(X,\cG)$.  Note that
$\varepsilon(X, \{\emptyset, \Omega\}) = \varepsilon(X)$
and $\bE[X \, | \, \{\emptyset, \Omega\}] = \bE[X]$.

More generally, suppose that $\cG \subseteq \cF$ is an arbitrary 
sub-$\sigma$-field and there is an associated regular conditional
probability $\bP_\cG(\omega', d \omega'')$ (such
a regular conditional probability certainly
exists if $\cG$ is finitely generated).  In this case,
we expect that $\bE[X \, | \, \cG](\omega')$ should be the
expectation of $X$ with respect to the probability measure
$\bP_\cG(\omega', \cdot)$.  It is easy to see that if
we let $\varepsilon(X, \cG)$ be the $\cG$-measurable random
variable such that $\varepsilon(X, \cG)(\omega')$ is
the infimum over $c \in \bK$ of the essential supremum
of $|X-c|$ with respect to $\bP_\cG(\omega', \cdot)$, then
this definition of $\varepsilon(X,\cG)$ subsumes our previous
one for the finitely generated case
and our putative definition of $\bE[X \, | \, \cG]$
coincides with the set of $\cG$-measurable random variables
$Y$ such that $|X-Y| \le \varepsilon(X,\cG)$, thereby
also extending the definition for the finitely
generated case.

We therefore see that the key to giving a satisfactory
general definition of $\bE[X \, | \, \cG]$ for an
arbitrary sub-$\sigma$-field $\cG \subseteq \cF$ is to
find a suitable general definition of $\varepsilon(X,\cG)$.
We tackle this problem in the next three sections.

\section{Conditional essential supremum}
\label{S:cond_ess_sup}

\begin{Def}
Given a non-negative real-valued random variable $S$ and
a sub-$\sigma$-field $\cG \subseteq \cF$, put
\[
\ess \sup\{S \, | \, \cG\} 
= \sup_{p \ge 1} \bE[S^p \, | \, \cG]^{\frac{1}{p}}
= \lim_{p \rightarrow \infty} \bE[S^p \, | \, \cG]^{\frac{1}{p}}.
\]
\end{Def}

\begin{Lemma}  
\label{L:cond_ess_sup_doms}
\begin{itemize}
\item[(i)]
Suppose that $S$ is a non-negative real-valued random variable
and $\cG$ is a sub-$\sigma$-field of $\cF$.  
Then $S \leq \ess\sup\{S\, | \, \cG\}$.
\item[(ii)]
Suppose that $S$ 
and $\cG$ are as in (i) and
$T$ is $\cG$-measurable real-valued
random variable with $S \le T$.
Then $\ess\sup\{S\, | \, \cG\} \leq T$.
\item[(iii)]
Suppose that $S'$ and $S''$ are 
non-negative real-valued random variables 
and $\cG$ is a sub-$\sigma$-fields of $\cF$.
Then 
\[
\ess\sup\{S' \vee S''\, | \, \cG\} = \ess\sup\{S'\, | \, \cG\} \vee \ess\sup\{S''\, | \, \cG\}.
\] 
\end{itemize}
\end{Lemma}
  
\begin{proof}
For part (i), we show by separate arguments that the result holds on the events
$\{\ess\sup\{S\, | \, \cG\}=0\}$ and 
$\{\ess\sup\{S\, | \, \cG\}>0\}$.

First consider what happens on the event
$\{\ess\sup\{S\, | \, \cG\}=0\}$.
By definition $\bE[S\, | \, \cG] \leq \ess\sup\{S\, | \, \cG\}$.  
Hence
\[
\begin{split}
\bE[S \bone\{\ess\sup\{S\, | \, \cG\}=0\}] 
& \leq \bE[S \bone\{\bE[S\, | \, \cG]=0\}] \\
& = \bE[\bE[S \bone\{\bE[S\, | \, \cG]=0\}\, | \, \cG]] \\
& =\bE[\bone\{\bE[S\, | \, \cG]=0\}\bE[S\, | \, \cG]] =0.\\
\end{split}
\]
Thus $\{\ess\sup\{S\, | \, \cG\}=0\} \subseteq \{S=0\}$, and  $S \leq \ess\sup\{S\, | \, \cG\}$ on the event
$\{\ess\sup\{S\, | \, \cG\}=0\}$.  
		
Now consider what happens on the event
$\{\ess\sup\{S\, | \, \cG\}>0\}$.
Take $\alpha > 1$.  Observe for $p \ge 1$ that 
\[
\begin{split}
\bE[S^p\, | \, \cG] 
& \geq \bE[S^p \bone\{S^p \geq \alpha^p\bE[S^p\, | \, \cG]\}\, | \, \cG] \\ 
& \geq \bE[\alpha^p\bE[S^p\, | \, \cG]\bone\{S^p\geq\alpha^p\bE[S^p\, | \, \cG]\}\, | \, \cG] \\
& = \alpha^p\bE[S^p\, | \, \cG] \; \bP\{S^p \geq \alpha^p\bE[S^p\, | \, \cG]\, | \, \cG\}. \\
\end{split}
\]
Hence, for each $p \ge 1$,
\[
\bP\{S \geq \alpha \:\ess\sup\{S\, | \, \cG\}\, | \, \cG\}
\le
\bP\{S \geq \alpha \: \bE[S^p\, | \, \cG]^{\frac 1 p}\, | \, \cG\} 
\leq \frac 1 {\alpha^p}
\]
on the event $\{\bE[S^p\, | \, \cG]>0\}$.  

Since $\{\ess\sup\{S\, | \, \cG\}>0\} \subseteq \bigcup_p \bigcap_{q \geq p} \{\bE[S^q\, | \, \cG]>0\}$, we see that  
$\bP\{S \geq \alpha \: \ess\sup\{S\, | \, \cG\}\, | \, \cG\}=0$
on the event on $\{\ess\sup\{S\, | \, \cG\}>0\}$.
As this holds for all $\alpha >1$, we conclude that 
$S \leq \ess\sup\{S\, | \, \cG\}$ on the event
$\{\ess\sup\{S\, | \, \cG\}>0\}$,
and this completes the proof of part (i).

Part (ii) is immediate from the definition. 

Now consider part (iii). We have from part (i) that 
$S' \le \ess\sup\{S'\, | \, \cG\}$ and $S'' \le  \ess\sup\{S''\, | \, \cG\}$.
Thus $S' \vee S'' \le \ess\sup\{S'\, | \, \cG\} \vee \ess\sup\{S''\, | \, \cG\}$
and hence 
\[
\ess\sup\{S' \vee S''\, | \, \cG\} \le \ess\sup\{S'\, | \, \cG\} \vee \ess\sup\{S''\, | \, \cG\}
\]
by part (ii).  On the other hand, because $S' \le S' \vee S''$ and $S'' \le S' \vee S''$, it follows
that $\ess\sup\{S'\, | \, \cG\} \le \ess\sup\{S' \vee S''\, | \, \cG\}$ and 
$\ess\sup\{S''\, | \, \cG\} \le \ess\sup\{S' \vee S''\, | \, \cG\}$. Therefore
\[
\ess\sup\{S'\, | \, \cG\} \vee \ess\sup\{S''\, | \, \cG\} \le \ess\sup\{S' \vee S''\, | \, \cG\}.
\]
\end{proof}

\begin{Cor}
\label{C:cond_reduces_Linfty}
Suppose that $S$ is a non-negative real-valued random variable
and $\cG \subseteq \cH$ are sub-$\sigma$-fields of $\cF$.
Then  
$\ess\sup\{S \, | \, \cH\} \leq \ess\sup\{S\, | \, \cG\}$.
\end{Cor}  

\begin{proof}
From Lemma \ref{L:cond_ess_sup_doms}(i),
$S \le \ess\sup\{S \, | \, \cG\}$.
Applying Lemma \ref{L:cond_ess_sup_doms}(ii)
with $\cG$ replaced by $\cH$ and $T = \ess\sup\{S \, | \, \cG\}$
gives the result.
\end{proof}

Let $\{\cF_n\}_{n=0}^\infty$ be a filtration (that is, a
non-decreasing sequence of sub-$\sigma$-fields of $\cF$).  
Recall that a random variable $T$ with values in 
$\{0,1,2,\ldots\}$ is a stopping time for the filtration
if $\{T=n\} \in \cF_n$ for all $n$.  Recall also
that if $T$ is a stopping time, then the associated
$\sigma$-field $\cF_T$ is the collection of 
events $A$ such that 
$A \cap \{T=n\} \in \cF_n$ for all $n$.

\begin{Lemma}
\label{L:ess_sup_stopping}
Suppose that $S$ is a non-negative real-valued random variable,
$\{\cF_n\}_{n=0}^\infty$ is a
filtration of sub-$\sigma$-fields of $\cF$, and $T$
is a stopping time.  Then 
\[
\begin{split}
& \ess\sup\{S \bone\{T=n\} \, | \, \cF_T\}
= \bone\{T=n\} \, \ess\sup\{S  \, | \, \cF_T\} \\
& \quad = \bone\{T=n\} \, \ess\sup\{S  \, | \, \cF_n\}
= \ess\sup\{S \bone\{T=n\} \, | \, \cF_n\} \\
\end{split}
\]
for all $n$.
\end{Lemma}

\begin{proof}
This follows immediately from the definition 
of the conditional essential supremum and the
fact that if 
$U$ is a non-negative real-valued random variable,
then
$\ess\sup\{U | \cF_T\} = \ess \sup\{U | \cF_n\}$ on the event $\{T=n\}$ (see, for example, Proposition II-1-3 of 
\cite{MR0402915}).
\end{proof}

\section{Conditional $L^\infty$ norm}
\label{S:cond_Linfty}

\begin{Def}
Given $X \in L^\infty$ and a sub-$\sigma$-field $\cG \subseteq \cF$,
put
\[
\|X\|_\cG := \ess \sup\{|X| \, | \, \cG\}.
\]
\end{Def}

\begin{Not}
Given $A \in \cF$, the $\bK$-valued
random variable $\bone_A$
is given by
\[
\bone_A(\omega) = 
\begin{cases}
1_{\bK},& \text{if $\omega \in A$},\\
0_{\bK},& \text{otherwise},
\end{cases}
\]
where $1_{\bK}$ and $0_{\bK}$ are, respectively, the multiplicative
and additive identity elements of $\bK$.  We 
continue to use this same notation to also
denote the analogously defined real-valued indicator random
variable, but this should cause no confusion as the
meaning will be clear from the context.
\end{Not}

\begin{Lemma}
\label{L:props_cond_Linfty}
Fix a sub-$\sigma$-field $\cG \subseteq \cF$.
\begin{itemize}
\item[(i)]
If $W \in L^\infty(\cG)$ and $X \in L^\infty$, 
then $\|W X\|_{\cG} = |W| \, \|X\|_{\cG}$.
\item[(ii)]
If $X, Y \in L^\infty$ are such that $\bP(\{X \ne Y\} \cap A) = 0$
for some $A \in \cG$, then $\bP(\{\|X\|_{\cG} \ne \|Y\|_{\cG}\} \cap A) = 0$.
\item[(iii)]
If $X_1,X_2,\ldots \in L^\infty$ and $A_1,A_2,\ldots \in \cG$ are pairwise disjoint, then 
\[
\left\|\sum_i X_i \bone_{A_i}\right\|_\cG = \sum_i \bone_{A_i} \|X_i\|_\cG.
\]
\item[(iv)]
If $X, Y \in L^\infty$, then
\[
\Vert X+Y \VertG \leq \Vert X \VertG \vee \Vert Y \VertG.
\]\end{itemize}
\end{Lemma}

\begin{proof}
Part (i) follows immediately from the definition. Part (ii) follows
from part (i): since $X \bone_A = Y \bone_A$ by assumption,
\[
\bone_A \|X\|_{\cG} = \|X \bone_A\|_{\cG} 
= \|Y \bone_A\|_{\cG} =\bone_A \|Y\|_{\cG}.
\]
Part (iii)
follows from parts (i) and (ii): for any of the events $A_j$,
\[
\begin{split}
\bone_{A_j} \sum_i \bone_{A_i} \|X_i\|_\cG
& =
\bone_{A_j} \|X_j\|_\cG
=
\|\bone_{A_j} X_j\|_\cG \\
& =
\left\|\bone_{A_j} \sum_i \bone_{A_i} X_i\right\|_\cG
=
\bone_{A_j} \left\|\sum_i \bone_{A_i} X_i\right\|_\cG, \\
\end{split}
\]
and, similarly, 
$\sum_i \bone_{A_i} \|X_i\|_\cG =\left\|\sum_i \bone_{A_i} X_i\right\|_\cG$
on $\Omega \setminus (\bigcup_i A_i)$.

Part (iv) is an immediate consequence of Lemma \ref{L:cond_ess_sup_doms}(iii).
However, there is also the following alternative, more elementary proof.
Note first that $\Vert X^r \VertG = \Vert X \VertG^r$ for any
$r > 0$ because
\[
\lim_{p\rightarrow\infty} \bE[|X|^{rp} \, | \, \cG]^{\frac 1 p} 
= \lim_{q\rightarrow\infty} \bE[|X|^q \, | \, \cG]^{\frac r q}
=(\lim_{q\rightarrow\infty} \bE[|X|^q \, | \, \cG]^{\frac 1 q})^r.
\]
Thus, from Jensen's inequality and  the observation that 
$(x+y)^s \leq (x^s+y^s)$ for $0 \leq s \leq 1$,
\[
\begin{split}
\Vert X+Y \VertG 
& = 
\lim_{p\rightarrow\infty} \bE[|X+Y|^p \, | \, \cG]^{\frac 1 p} \\
& \leq  
\lim_{p\rightarrow\infty} \bE[|X|^p \vee |Y|^p\, | \, \cG]^{\frac 1 p} \\
& = 
\lim_{p\rightarrow\infty} \bE[\lim_{r\rightarrow\infty} (|X|^{pr}+|Y|^{pr})^{\frac 1 r}\, | \, \cG )^{\frac 1 p} \\
& \leq \lim_{p,r\rightarrow\infty}(\bE[|X|^{pr}\, | \, \cG]+\bE[|Y|^{pr}\, | \, \cG])^{\frac 1 {pr}} \\
& \leq \lim_{p,r\rightarrow\infty}(\bE[|X|^{rp}\, | \, \cG]^{\frac 1 p}+\bE[|Y|^{rp}\, | \, \cG]^{\frac 1 p})^{\frac 1 r}.\\
& = \lim_{r\rightarrow\infty}(\Vert X \VertG^r+\Vert Y \VertG^r)^{\frac 1 r} \\
& = \Vert X \VertG \vee \Vert Y \VertG.\\
\end{split}
\]
\end{proof}

The following result is immediate
from Corollary \ref{C:cond_reduces_Linfty}.

\begin{Lemma}
\label{L:cond_reduces_norm}	
Suppose that $X \in L^\infty$  and 
$\cG \subseteq \cH$ are sub-$\sigma$-fields of $\cF$.
Then $\|X\|_\cH \leq \|X\|_\cG$. 
\end{Lemma}

The following result is immediate
from Lemma \ref{L:ess_sup_stopping}.

\begin{Lemma}
\label{L:Linfty_stopping}
Suppose that $X \in L^\infty$,
$\{\cF_n\}_{n=0}^\infty$ is a
filtration of sub-$\sigma$-fields of $\cF$, and $T$
is a stopping time.  Then 
\[
\begin{split}
& \|X \bone\{T=n\}\|_{\cF_T}
= \bone\{T=n\} \, \|X\|_{\cF_T} \\
& \quad = \bone\{T=n\} \, \|X\|_{\cF_n}
= \|X \bone\{T=n\}\|_{\cF_n} \\
\end{split}
\]
for all $n$.
\end{Lemma}

\section{Construction of Conditional Expectation}
\label{S:constr_cond_exp}

	\begin{Def}
Given $X \in L^\infty$ and a sub-$\sigma$-field $\cG \subseteq \cF$,
set
\[
		\bE[X \, | \, \cG] := \{Y \in L^\infty(\cG): \Vert X-Y \VertG \leq \Vert X-Z \VertG \text{ for all } Z \in L^\infty(\cG) \}.
\]
	\end{Def}
	
\begin{Rem}
Before showing that $\bE[X \, | \, \cG]$ is non-empty, we comment
on a slight subtlety in the definition.  One way of thinking of
our definition of $\bE[X]$ as the set of $c \in \bK$
for which $\|X - c\|_\infty$ is minimal, 
is that $\bE[X]$ is the set of projections
of $X$ onto $\bK \equiv L^\infty(\{\emptyset, \Omega\})$.  
A possible definition
of $\bE[X \, | \, \cG]$ might therefore be the analogous set of
projections of $X$ onto $L^\infty(\cG)$, that is, the set
of $Y \in L^\infty(\cG)$ that minimize $\|X-Y\|_\infty$.
This definition is {\bf not} equivalent to ours.  For example,
suppose that $\Omega$ consists of the three points 
$\{\alpha, \beta, \gamma\}$, 
$\cF$  consists of all subsets of $\Omega$, $\bP$ assigns positive
mass to each point of $\Omega$, $\cG = \sigma\{\{\alpha, \beta\}, \{\gamma\}\}$, and $X$ is given by $X(\alpha) = 1_\bK$, $X(\beta) = 0_\bK$,
and $X(\gamma) = 0_\bK$.  Consider $Y \in L^\infty(\cG)$,
so that $Y(\alpha) = Y(\beta) = c$ and
$Y(\gamma) = d$ for some $c,d \in \bK$.  In order that
$Y \in \bE[X \, | \, \cG]$
according to our definition, $c$ and $d$ must be chosen to minimize both
$|1_\bK - c| \vee |0_\bK - c|$ and $|0_\bK - d|$.  By the strong triangle
inequality, $|1_\bK - c| \vee |0_\bK - c|$ is minimized by any
$c$ with $|c| \le 1$, with the corresponding minimal value being $1$.
Of course, $|0_\bK - d|$ is minimized by the unique value $d=0_K$.
On the other hand, in order that $Y$ is a projection of $X$
onto $L^\infty(\cG)$, the points $c$ and $d$ must be chosen to minimize
$|1_\bK - c| \vee |0_\bK - c| \vee |0_\bK - d|$, and this is
accomplished as long as $|c| \le 1$ and $|d| \le 1$.  We don't
belabor the point in what follows, but several of the natural counterparts
of standard results for classical conditional expectation
that we show hold for our definition fail to hold for the
``projection'' definition.
\end{Rem}

The following lemma is used below to show that 
$\bE[X \, | \, \cG]$ is non-empty.

	\begin{Lemma}
	\label{L:projective}
	Suppose that $X \in L^\infty$ is not $0_\bK$ almost surely, and $\cG$
	is a sub-$\sigma$-field of $\cF$.  Set $q^{-N} = \|X\|_\infty$.
	Then there exist disjoint events $A_0, A_1, \ldots \in \cG$ 
	and random variables
	$Y_0, Y_1, \ldots \in L^\infty(\cG)$ with the 
	following  properties:
	\begin{itemize}
	\item[(1)] On the event $A_n$, $\|X - Z\|_\cG \ge q^{-(N+n)}$
	for every  $Z \in L^\infty(\cG)$. 
	\item[(2)] On the event $A_n$, 
	$\|X - Y_n\|_\cG = q^{-(N+n)}$.  and 
	\item[(3)] On the event $\Omega \setminus  \bigcup_{k=1}^n A_k$,
	$\|X - Y_n\|_\cG \le q^{-(N+n+1)}$
  \item[(4)] On the event $\bigcup_{k=1}^n A_k$, $Y_p = Y_n$ for 
  any $p > n$.
  \item[(5)] The event $\bigcup_{k=1}^\infty A_k$ has probability one. 
	\end{itemize}
	\end{Lemma}

		\begin{proof}
		Suppose without loss of generality that $\|X\|_\infty = 1$,
		so that $N=0$.
		Set 
		$\mathbf{Z_0} := \{Z \in L^\infty(\cG) : \|X-Z\|_\infty \le 1\}$.
		Note that the constant $0$
		belongs to $\mathbf{Z_0}$ and so this set is non-empty.
		Put
		$\delta_0 := \inf_{Z \in \mathbf{Z_0}}
		\bP\{\|X - Z\|_\cG = 1\}$.  

		Choose $Z_{0,1}, Z_{0,2}, \ldots \in \mathbf{Z_0}$ with
		\[
		\lim_{m \rightarrow \infty}
		\bP\{\|X - Z_{0,m}\|_\cG = 1\}
		=
		\delta_0.
		\]
		Define $Z_{0,1}', Z_{0,2}', \ldots$ inductively by setting 
		$Z_{0,1}' := Z_{0,1}$ and
		\[
		Z_{0, m+1}'(\omega) := 
		\begin{cases}
		Z_{0,m}'(\omega),& \text{if $\|X - Z_{0,m}'\|_\cG(\omega) \le \|X - Z_{0, m+1}\|_\cG(\omega)$}, \\
		Z_{0, m+1}(\omega),& \text{if $\|X - Z_{0,m}'\|_\cG(\omega) > \|X - Z_{0, m+1}\|_\cG(\omega)$}.
		\end{cases}
		\]

		Note that the events $B_{0,m} := \{\|X - Z_{0,m}'\|_\cG = 1\}$ are decreasing and the $B_{0,m}$ are contained in the event 
$\{\|X - Z_{0,m}\|_\cG = 1\}$.  
Hence the event $A_0 := \lim_{m \rightarrow \infty} B_{0,m} = \bigcap_{m=1}^\infty B_{0,m}$ has probability $\delta_0$.

Define $Y_0$ by
		\[
		Y_0(\omega) := 
		\begin{cases}
		Z_{m,1}'(\omega), & \text{if $\omega \in (\Omega \setminus  B_{0,1}) \cup A_0$}, \\
		Z_{0,m}'(\omega), & \text{if $\omega \in (\Omega \setminus  B_{0,m})
		\setminus  (\Omega \setminus  B_{0,m-1}), 
		\; m \ge 2$}.
		\end{cases}
		\]

		It is clear that
		$\|X - Y_0\|_\cG = 1$ on the event $A_0$
		and $\|X - Y_0\|_\cG \le q^{-1}$ on the event
		$\Omega \setminus  A_0$.
		Moreover, if there existed $V \in L^\infty(\cG)$
		with
		\[
		\bP(\{\|X - V\|_\cG \le q^{-1}\} \cap A_0) > 0,
		\]
		then we would have the contradiction that
		$W \in \mathbf{Z}_0$ defined by
		\[
		W(\omega)
		= 
		\begin{cases}
		Y_0(\omega), & \text{if $\|X - Y_0\|_\cG(\omega) \le \|X - V\|_\cG(\omega)$,} \\
		V(\omega), & \text{if $\|X - Y_0\|_\cG(\omega) > \|X - V\|_\cG(\omega)$}
		\end{cases}
		\]
		would satisfy $\bP\{\|X - W\|_\cG=1\} < \delta_0$.

		Now suppose that $A_0, \ldots A_{n-1}$ and
		$Y_0, \ldots, Y_{n-1}$ have been constructed with the
		requisite properties. If
		$\bP(\Omega \setminus \bigcup_{k=1}^{n-1})=0$, then
		take $A_n = \emptyset$ and $Y_n = Y_{n-1}$ (recall that
		we are interpreting all equalities and inequalities as
		holding $\bP$-a.s.)  Otherwise,
		set 
\[
\begin{split}
\mathbf{Z_n} := 
\biggl\{Z \in L^\infty(\cG) : 
& \text{$Z = Y_{n-1}$ on $\bigcup_{k=1}^{n-1} A_k$} \\
& \text{and $|X-Z| \le q^{-n}$ on 
		$\Omega \setminus \bigcup_{k=1}^{n-1} A_k$}\biggl\}.\\
\end{split}
\]  
		Note
		that $Y_{n-1}$ belongs to $\mathbf{Z_n}$.
		Put
		$\delta_n := \inf_{Z \in \mathbf{Z_n}}
		\bP\{\|X - Z\|_\cG = q^{-n}\}$. 
		An argument very similar to the above with 
		$\mathbf{Z_n}$ and $\delta_n$ replacing 
		$\mathbf{Z_0}$ and $\delta_0$ establishes the
		existence of $A_n$ and $Y_n$ with the desired properties.
		\end{proof}
		
			\begin{Theo}
			Given $X \in L^\infty$ and 
a sub-$\sigma$-algebra $\cG \subseteq \cF$, the conditional
expectation $\bE[X \, | \, \cG]$ is nonempty.
	\end{Theo}
	
\begin{proof} If $X$ is $0_\bK$ almost surely, then
$\bE[X \, | \, \cG] = \{0_\bK\}$.  Otherwise, let
$A_0, A_1, \ldots \in \cG$ 
and
$Y_0, Y_1, \ldots \in L^\infty(\cG)$
be as in Lemma \ref{L:projective}. Then $Y$ defined by
$Y(\omega) = Y_n(\omega)$ for $\omega \in A_n$
belongs to $\bE[X \, | \, \cG]$.
\end{proof}

\section{Elementary Properties of Conditional Expectation}
\label{S:prop_cond_exp}
	
	\begin{Prop}
\label{P:multiplication_addition}
Fix a sub-$\sigma$-field $\cG \subseteq \cF$.
\begin{itemize}
\item[(i)]
Suppose that $X \in L^\infty(\cG)$ and $Y \in L^\infty$.
Then
\[
		\bE[XY\, | \, \cG] = X \, \bE[Y\, | \, \cG].
\]
and
\[
		\bE[X+Y\, | \, \cG] = X + \bE[Y\, | \, \cG].
\]
\item[(ii)]
If $X, Y \in L^\infty$ are such that $\bP(\{X \ne Y\} \cap A) = 0$
for some $A \in \cG$, then 
$\bone_A \bE[X \, | \, \cG] = \bone_A \bE[ Y \, | \, \cG]$.
\item[(iii)]
If $X_1,X_2,\ldots \in L^\infty$ and $A_1,A_2,\ldots \in \cG$ are pairwise disjoint, then 
\[
\bE\left[\sum_i X_i \bone_{A_i} \, | \, \cG\right] = \sum_i \bone_{A_i} \bE[X_i \, | \, \cG].
\]
\end{itemize}
\end{Prop}
	
		\begin{proof}
Consider part (i).		We first show the inclusion
$\bE[XY\, | \, \cG] \subseteq X\bE[Y\, | \, \cG]$.  
		
		Consider $Z \in \bE[XY\, | \, \cG]$.  
Choose some $V \in \bE[Y\, | \, \cG]$ and set 
$W = (Z/X) \bone\{X \ne 0\} + V \bone\{X = 0\} \in L^\infty(\cG)$. Note
that $\bP\{Z \ne 0, \, X=0\} = 0$ and hence $X W = Z$, 
because otherwise we would have the contradiction
		$\Vert XY-Z \bone\{X \ne 0\} \VertG \le \Vert XY-Z \VertG$
		and
		$\bP\{\Vert XY-Z \bone\{X \ne 0\} \VertG < \Vert XY-Z \VertG\}>0$ by Lemma \ref{L:props_cond_Linfty}(ii).
		
We need to show that $W \in \bE[Y\, | \, \cG]$.
Consider $U \in L^\infty(\cG)$.
By Lemma \ref{L:props_cond_Linfty}(ii) and the assumption
that $V \in \bE[Y\, | \, \cG]$,
\[
\Vert Y - W \VertG 
= \Vert Y - V \VertG 
\le \Vert Y - U \VertG
\] on the event $\{X=0\}$. Also,
$\Vert XY-Z \VertG \le \Vert XY-XU \VertG$ 
by the assumption that $Z \in \bE[XY\, | \, \cG]$,
and so, by Lemma \ref{L:props_cond_Linfty}(i)$+$(ii)
\[
\begin{split}
\Vert Y - W \VertG & = \Vert Y - Z/X \VertG
= |X|^{-1} \Vert XY-Z \VertG \\
& \le |X|^{-1} \Vert XY-XU \VertG
= \Vert Y-U \VertG \\
\end{split}
\]
on the event $\{X \ne 0\}$.  Thus
$\Vert Y - W \VertG \le \Vert Y - U \VertG$
for any $U \in L^\infty(\cG)$ and hence
$W \in \bE[Y\, | \, \cG]$.

We now show the converse inclusion 
$X\bE[Y\, | \, \cG] \subseteq \bE[XY\, | \, \cG]$.  
		
Choose $W \in \bE[Y\, | \, \cG]$.  We need to show
that $X W \in \bE[XY\, | \, \cG]$.  Consider $U \in L^\infty(\cG)$.   Put $V = (U/X) \bone\{X \ne 0\}$. 
We have 
$\Vert Y - W \VertG
\le \Vert Y -  V \VertG$ by the assumption that
$W \in \bE[Y\, | \, \cG]$.
From Lemma \ref{L:props_cond_Linfty}(i)$+$(ii),
\[
\begin{split}
\Vert XY - XW \VertG 
& = |X| \Vert Y - W \VertG
\le
|X| \Vert Y -  V \VertG
= \Vert XY - XV \VertG \\
& = \Vert XY - U \VertG \bone\{X \ne 0\}
\le \Vert XY - U \VertG, \\
\end{split}
\]
as required.	

The proof of the claim $\bE[X+Y\, | \, \cG] = X + \bE[Y\, | \, \cG]$
is similar but easier, so we omit it.

Parts (ii) and	(iii) follow straightforwardly from parts (ii) and (iii)
of Lemma \ref{L:props_cond_Linfty}.
\end{proof}

	\begin{Prop}
\label{P:independence_cond_exp}
	Let $\cG$ be a sub-$\sigma$-algebra of $\cF$.
Suppose that $X \in L^\infty$  is independent of $\cG$. Then 
$\bE[X \, | \, \cG]$ is the set of random variables
$Y \in L^\infty(\cG)$ that take values in $\bE[X]$.
\end{Prop}

		\begin{proof}
Observe for any $Z \in L^\infty(\cG)$, 
that, by the assumption of independence of $X$ from $\cG$, 
\[
\begin{split}
\|X-Z\|_\cG(\omega) 
& = \sup_p \left(\bE[|X-Z|^p\, | \, \cG](\omega)\right)^{\frac 1 p} \\
& = \sup_p \left(\int |x-Z(\omega)|^p \;\bP\{X \in dx\}\right)^{\frac 1 p} \\
& = \sup\{|x-Z(\omega)|:x \in \supp X\} \\
& 
\begin{cases} 
= \varepsilon(X),& \text{ if $Z(\omega) \in \bE[X]$}, \\
> \varepsilon(X),& \text{otherwise},
\end{cases} \\
\end{split}
\]
and the result follows.
		\end{proof}

\section{Conditional spread and the tower property}
\label{S:tower_prop}

\begin{Def}
Given $X \in L^\infty$ and a sub-$\sigma$-field $\cG$ of $\cF$, let
 $\varepsilon(X,\cG)$ denote the common value of $\|X-Y\|_\cG$ for 
 $Y \in \bE[X \, | \, \cG]$.
 \end{Def} 
 
\begin{Lemma}
\label{L:cond_smaller_spread}
If $X \in L^\infty$ and a $\cG \subseteq \cH$
are sub-$\sigma$-fields of $\cF$, then
$\varepsilon(X,\cH) \le \varepsilon(X,\cG)$.
\end{Lemma}

\begin{proof}
Suppose that $V \in \bE[X \, | \, \cG]$ and $W \in \bE[X \, | \, \cH]$.
From Lemma \ref{L:cond_reduces_norm},
\[
\varepsilon(X, \cH) 
= \Vert X - W\Vert_\cH 
\le \Vert X - V\Vert_\cH
\le \Vert X - V\Vert_\cG
= \varepsilon(X, \cG).
\]
\end{proof}

	\begin{Lemma}
\label{L:alternative_def_cond_exp}
A random variable $Y$ belongs to $\bE[X \, | \, \cG]$ 
if and only if $Y \in L^\infty(\cG)$  and $|X-Y| \leq \varepsilon(X,\cG)$.
	\end{Lemma}
	
		\begin{proof}
Suppose $Y$ is in $\bE[X \, | \, \cG]$. By definition, $Y \in L^\infty(\cG)$.  
By Lemma Lemma \ref{L:cond_reduces_norm},
$|X-Y| = \|X-Y\|_\cF \leq \|X-Y\|_\cG = \varepsilon(X, \cG)$.  

The converse is immediate from Lemma \ref{L:cond_ess_sup_doms}(ii).
		\end{proof}

	\begin{Lemma} 
\label{L:tower_spread}
Suppose that	$X \in L^\infty$,  $\cG \subseteq \cH$
are sub-$\sigma$-fields of $\cF$,
and $Y \in \bE[X \, | \, \cH]$.  Then $\varepsilon(Y,\cG) \leq \varepsilon(X,\cG)$.
	\end{Lemma} 
	 
		\begin{proof}
Consider $Z \in \bE[X \, | \, \cG]$.  By 
Lemma \ref{L:alternative_def_cond_exp} and 
Lemma \ref{L:cond_smaller_spread}
		\[
		|Y-Z| \leq |X-Y| \vee |X-Z| \leq \varepsilon(X,\cH) \vee \varepsilon(X,\cG) = \varepsilon(X,\cG).
		\] 
By Lemma \ref{L:cond_ess_sup_doms}(ii),
$\varepsilon(Y,\cG) \le \|Y-Z\|_\cG \le \varepsilon(X,\cG)$.
\end{proof}

\begin{Theo}
\label{T:tower_prop}
Suppose that	$X \in L^\infty$ and  $\cG \subseteq \cH$
are sub-$\sigma$-fields of $\cF$.
If $Y \in \bE[X \, | \, \cH]$ and
	$Z \in \bE[Y\, | \, \cG]$, then $Z \in \bE[X \, | \, \cG]$.
	\end{Theo}

\begin{proof}
By Lemma \ref{L:alternative_def_cond_exp},  
Lemma \ref{L:tower_spread}, and Lemma \ref{L:cond_smaller_spread},
	\[
	|X-Z| \leq |X-Y| \vee |Y-Z| \leq \varepsilon(X,\cH) \vee \varepsilon(Y,\cG) \le  \varepsilon(X,\cG).
	\]
Thus $Z$ is in $\bE[X \, | \, \cG]$
by another application of Lemma \ref{L:alternative_def_cond_exp}.
\end{proof} 
		
\section{Continuity of conditional expectation}
\label{S:conty_cond_exp}

\begin{Def}
Define the Hausdorff distance between two subsets $A$ and $B$
of $L^\infty$ to be
\[
D_H(A,B) :=
\sup_{X \in A} \inf_{Y \in B} \|X-Y\|_\infty 
\vee
\sup_{Y \in B} \inf_{X \in A} \|Y-X\|_\infty.
\]
\end{Def}

\begin{Lemma}
\label{L:Hausdorff_Minkowski_sum}
Suppose that $A,B,C$ are subsets of $L^\infty$.  Then
\[
D_H(A+C,B+C) \leq D_H(A,B).
\]
\end{Lemma}

\begin{proof}
Suppose that $D_H(A,B) < \delta$ for some $\delta \ge 0$.
By definition,
 for every $X \in A$ there is a $Y \in B$ with 
$\|X-Y\|_\infty < \delta$, and similarly with the roles of $A$ and $B$ reversed.  
If $U \in A+C$, then $U = X+W$ for some $X \in A$ and $W \in C$.  We know there is $Y \in B$ such that 
$\|X-Y\|_\infty < \delta$.  
Note that $V := Y + W \in B + C$ and 
$\|U-V\|_\infty = \|X-Y\|_\infty < \delta$.
A similar argument with the roles of $A$ and $B$ reversed shows that $D_H(A+C,B+C) < \delta$.
\end{proof}

\begin{Theo}
\label{T:conty_cond_exp}
Suppose that $X,Y \in L^\infty$ and $\cG$
is a sub-$\sigma$-field of $\cF$. Then
$D_H(\bE[X \, | \, \cG], \bE[Y\, | \, \cG])
\le \|X - Y\|_\infty$ .
\end{Theo}
		
\begin{proof}
Choose $U \in\bE[X \, | \, \cG]$  and
$V \in \bE[Y\, | \, \cG]$.  From Lemma \ref{L:props_cond_Linfty}(iv),
		\[
		\varepsilon(Y,\cG) \le \|Y-U\|_\cG \le \|X-U\|_\cG \vee \|X-Y\|_\cG = \varepsilon(X,\cG) \vee \|X-Y\|_\cG
		\]
and 
		\[
		\varepsilon(X,\cG) \le \|X-V\|_\cG \le \|Y-V\|_\cG \vee \|X-Y\|_\cG = \varepsilon(Y,\cG) \vee \|X-Y\|_\cG.
		\]

It follows that 
$\varepsilon(X,\cG) = \varepsilon(Y,\cG)$
on the event	
$M := \{\|X-Y\|_\cG < \varepsilon(X,\cG) \vee \varepsilon(Y,\cG)\}$ and 	
\[
		\varepsilon(X,\cG) 
		= \|Y-U\|_\cG
		= \|X-U\|_\cG  
		= \varepsilon(X,\cG)
\]
		and
\[
		\varepsilon(Y,\cG) 
		= \|X-V\|_\cG 
		= \|Y-V\|_\cG  
		= \varepsilon(Y,\cG)
\]
on $M$.		
		
By Proposition \ref{P:multiplication_addition},
$U \bone_M \in \bE[Y \bone_M\, | \, \cG]=\bone_M \bE[Y\, | \, \cG]$ and $V \bone_M \in \bE[X \bone_M\, | \, \cG]=\bone_M \bE[X \, | \, \cG]$. 
Thus
 $\bone_M \bE[X \, | \, \cG] = \bone_M \bE[Y \, | \, \cG]$.
		
Furthermore, on the event
$N:= \{\|X-Y\|_\cG \geq \varepsilon(X,\cG) \vee \varepsilon(Y,\cG)\}$ 
		\[
		\begin{split}
		\|U-V\|_\infty 
		& \le \|U - X\|_\infty \vee \|X-Y\|_\infty \vee \|Y - V\|_\infty \\
		& \le \varepsilon(X,\cG) \vee \|X-Y\|_\infty \vee \varepsilon(Y,\cG) \\
		& \le \|X-Y\|_\infty, \\
		\end{split}
		\]
and so
$\|U \bone_N - V \bone_N\|_\infty \leq \|X \bone_N - Y \bone_N\|_\infty \leq \|X-Y\|_\infty$.  Therefore,
\[
D_H(\bone_N \bE[X \, | \, \cG], \bone_N \bE[Y \, | \, \cG]) \le \|X-Y\|_\infty.
\]
		
By Proposition \ref{P:multiplication_addition}(iii),
$\bE[X \, | \, \cG] = \bone_M \bE[X \, | \, \cG]+ \bone_N \bE[X \, | \, \cG]$, and similarly for $Y$.  The
result now follows from 
Lemma \ref{L:Hausdorff_Minkowski_sum}. 
\end{proof}

\section{Martingales}
\label{S:mart}

		\begin{Def}
Let $\{\cF_n\}_{n=0}^\infty$ be a
filtration of sub-$\sigma$-fields of $\cF$.
A sequence of random variables $\{X_n\}_{n=0}^\infty$
is a {\em martingale} if there exists $X \in L^\infty$ such that
$X_n \in \bE[X \, | \, \cF_n]$ for all $n$ (in particular, 
$X_n \in L^\infty(\cF_n)$).
  \end{Def}
  
\begin{Rem}
Note that our definition does not imply that 
$X_n \in \bE[X_{n+1} \, | \, \cF_n]$ for all $n$.
For example, suppose that $\cF_n := \{\emptyset, \Omega\}$
for all $n$ but $X$ is not almost surely constant, then
we obtain a martingale by taking $X_n$ to be any constant in the
ball $\bE[X]$, but we only have $X_n \in \bE[X_{n+1} \, | \, \cF_n]$
for all $n$ if $X_0=X_1=X_2= \ldots$.  
\end{Rem}

Many of the usual real-valued examples of martingales
have $\bK$-valued counterparts.

\begin{Exm}
Let $\{Y_n\}_{n=0}^\infty$ be a sequence
of independent random variables in $L^\infty$
with $0_\bK \in \bE[Y_n]$ for all $n$. 
Suppose that $\sum_{k=0}^\infty Y_k$
converges in $L^\infty$ (by the
strong triangle inequality and the
completeness of $L^\infty$, this is equivalent
to $\lim_{n \rightarrow \infty} \|Y_n\|_\infty = 0$).
Set $\cF_n := \sigma\{Y_0, Y_1, \ldots, Y_n\}$.  Put
$X_n := \sum_{k=0}^n Y_k$ and
$X_n := \sum_{k=0}^\infty Y_k$
It follows from
the second claim of 
Proposition \ref{P:multiplication_addition}(i) that
$X_n \in \bE[X \, | \, \cF_n]$ for all $n$ and hence
$\{X_n\}_{n=0}^\infty$ is a martingale.
\end{Exm}

\begin{Exm}
Let $\{Y_n\}_{n=0}^\infty$ be a sequence
of independent random variables in $L^\infty$
with $1_\bK \in \bE[Y_n]$ for all $n$. 
Suppose that $\prod_{k=0}^\infty Y_k$
converges in $L^\infty$ (by the
strong triangle inequality and the
completeness of $L^\infty$, this is equivalent
to $\lim_{n \rightarrow \infty} \|Y_n - 1_\bK\|_\infty = 0$).
Set $\cF_n := \sigma\{Y_0, Y_1, \ldots, Y_n\}$.
Put $X_n := \prod_{k=0}^n Y_k$ and $X := \prod_{k=0}^\infty Y_k$.  
It follows from
the first claim of 
Proposition \ref{P:multiplication_addition}(i) that
$X_n \in \bE[X \, | \, \cF_n]$ for all $n$ and hence
$\{X_n\}_{n=0}^\infty$ is a martingale.
\end{Exm}

\begin{Exm}
Let $\{Z_n\}_{n=0}^\infty$ be a discrete time
Markov chain with countable state space $E$
and transition matrix $P$.  Set
$\cF_n := \sigma\{Z_0, Z_1, \ldots, Z_n\}$.
Say that $f:E \rightarrow \bK$
is {\em harmonic} if $f$ is
bounded and for all $i \in E$ the
expectation of $f$ with respect to
the probability measure $P(i,\cdot)$ contains $f(i)$
(that is, if $f(i)$ is belongs to the smallest ball
containing the set $\{f(j) : P(i,j) > 0\}$).
Fix $N \in \{0,1,2,\ldots\}$.
Then $\{X_n\}_{n=0}^\infty := \{f(Z_{n \wedge N})\}_{n=0}^\infty$
is a martingale.
\end{Exm}

\section{Optional sampling theorem}
\label{S:opt_sampling}

		\begin{Theo}
		Let $\{\cF_n\}_{n=0}^\infty$ be a
filtration. Suppose that $X \in L^\infty$ and
$\{X_n\}_{n=0}^\infty$ is a martingale
with $X_n \in \bE[X \, | \, \cF_n]$ for all $n$.
If $T$ is a stopping time, then
$X_T \in \bE[X \, | \, \cF_T]$. 
\end{Theo} 
			
			\begin{proof}
It follows from Lemma \ref{L:Linfty_stopping} that
$\bone\{T=n\} \bE[X  \, | \, \cF_T] 
= \bone\{T=n\} \bE[X  \, | \, \cF_n]$ and hence,
by Proposition \ref{P:multiplication_addition}(iii),
\[
\begin{split}
\bE[X  \, | \, \cF_T]
& =
\bE\left[\sum_n X \bone\{T=n\} \, | \, \cF_T\right] \\
& =
\sum_n \bone\{T=n\} \bE[X  \, | \, \cF_T] \\
& =
\sum_n \bone\{T=n\} \bE[X  \, | \, \cF_n] \\
& \ni
\sum_n \bone\{T=n\} X_n \\
& =
X_T. \\
\end{split}
\]
\end{proof}

\section{Martingale convergence}
\label{S:mart_conv}

	\begin{Theo}

		Let $\{\cF_n\}_{n=0}^\infty$ be a
filtration. Suppose that $X \in L^\infty$ and
$\{X_n\}_{n=0}^\infty$ is a martingale
with $X_n \in \bE[X \, | \, \cF_n]$ for all $n$.
If $X$ is in the closure of $\bigcup_{n=1}^{\infty}L^\infty(\cF_n)$,
then $\lim_{n \rightarrow \infty} \|X_n - X\|_\infty = 0$
(in particular, $\{X_n\}_{n=0}^\infty$ converges to
$X$ almost surely).
\end{Theo}

\begin{proof}
		Since $X$ is in the closure of
$\bigcup_{n=1}^{\infty}L^\infty(\cF_n)$, for each $\varepsilon > 0$ there exists 
$Y \in L^\infty(\cF_N)$ for some $N$ such that 
$\|X-Y\|_\infty < \varepsilon$.  Because $\cF_N \subseteq \cF_n$
for $n > N$, $Y \in L^\infty(\cF_n)$ for $n \ge N$.

By Theorem \ref{T:conty_cond_exp},
$D_H(E[X \, | \, \cF_n], E[Y \, | \, \cF_n]) < \varepsilon$ for
$n \ge N$.  However, $E[Y \, | \, \cF_n]$
consists of the single point $Y$, and so
the Hausdorff distance is simply
$\sup\{\|W-Y\|_\infty : W \in E[X \, | \, \cF_n]\}$.
Thus
\[
\|X_n - X\|_\infty 
\le \|X_n - Y\|_\infty \vee \|Y - X \|_\infty 
< \varepsilon
\]
for $n \ge N$.
\end{proof}

\def\polhk#1{\setbox0=\hbox{#1}{\ooalign{\hidewidth
  \lower1.5ex\hbox{`}\hidewidth\crcr\unhbox0}}}
  \def\polhk#1{\setbox0=\hbox{#1}{\ooalign{\hidewidth
  \lower1.5ex\hbox{`}\hidewidth\crcr\unhbox0}}}
  \def\polhk#1{\setbox0=\hbox{#1}{\ooalign{\hidewidth
  \lower1.5ex\hbox{`}\hidewidth\crcr\unhbox0}}}
\providecommand{\bysame}{\leavevmode\hbox to3em{\hrulefill}\thinspace}
\providecommand{\MR}{\relax\ifhmode\unskip\space\fi MR }
\providecommand{\MRhref}[2]{%
  \href{http://www.ams.org/mathscinet-getitem?mr=#1}{#2}
}
\providecommand{\href}[2]{#2}

\end{document}